\documentclass[11pt]{article}
\newif\ifpdf
\ifx\pdfoutput\undefined
    \pdffalse       
\else
    \pdfoutput=1    
    \pdftrue
\fi

\usepackage{amssymb,amsmath,amsthm,url}

\ifpdf
    \usepackage[pdftex]{graphicx}
    \usepackage[colorlinks=false, pdfstartview=FitH, linkcolor=blue,
        citecolor=blue, urlcolor=blue]{hyperref}
    \usepackage{amsrefs}
\else
    \usepackage{graphicx}
    \usepackage{hyperref}
    \usepackage{amsrefs}
\fi

    \newcommand{\MathReview}[1]{}

    \newcommand{\N}{\mbox{$\mathbb N$}}

    \newcommand{\floor}[1]{{\left\lfloor #1 \right\rfloor}}
    
    \newcommand{\fp}[1]{{\{ #1 \}}}
    \newcommand{\tf}[1]{{[\![ #1 ]\!]}}
    \newcommand{\bigO}[1]{{\mathcal O\left( #1 \right)}}
    \newcommand{\tbigO}[1]{{\mathcal O\big( #1 \big)}}

    \newcommand{\al}{\alpha}
    
    \newcommand{\be}{\beta}

    \newtheorem{thm}{Theorem}[section]
    \newtheorem{idea}{Idea}
    \newtheorem{lem}[thm]{Lemma}
    \newtheorem{cor}[thm]{Corollary}

    \newenvironment{proofof}[1]{\medskip\noindent{\em Proof of #1.}}{\qed\medskip}

    \renewcommand{\section}[1]{\vspace{12pt} \addtocounter{section}{1} \noindent{\bf \thesection. #1} \,}
    \renewcommand{\subsection}[1]{\vspace{12pt} \addtocounter{subsection}{1} \noindent {\bf #1} \,}

\title{Almost Alternating Sums}
\author{
  Kevin O'Bryant,
  Bruce Reznick, and
  Monika Serbinowska}
\date{}

\begin{document}

\maketitle



\section{INTRODUCTION.}
The behavior of the sum
    $$S_N(\al) := \sum_{n=1}^N (-1)^\floor{n\al}$$
as $N\to\infty$ is not transparent (here $\floor{x}$ signifies the greatest integer not larger than $x$). The
random walk
    $ \sum_{n=1}^N w_n,$
where the $w_n$ are independent random variables taking the values $1$ and $-1$ with equal probability, is known
\cite{MathWorld.RandomWalk} to typically have absolute value around $c\sqrt{N}$ for an appropriate constant $c$
and large $N$. Knowing this and knowing also that for irrational $\al$ the sequence $\floor{n\al}$ is
``random-ish'' modulo 2, a natural guess is that $|S_N(\al)|$ is also around $\sqrt{N}$.

Contrary to this expectation, for almost all real numbers $\al$
    \begin{equation}\label{eq:almostall}
    \big|S_N(\al)\big| \le ( \log N )^2
    \end{equation}
for all large $N$. This is a corollary of a theorem of Khintchine, which we state precisely in section~2.

We devote the bulk of this article to two elementary proofs that
    \begin{equation}\label{eq.main}
    \big|S_N(\al)\big| \le \frac {\log N}{2\log (1 + \sqrt 2)} + 1
    \end{equation}
for all $N$ and an explicit countable set of $\al$, including $\sqrt{2}$ and $\sqrt{5}+1$. Our first proof is
entirely self-contained and is given in section~3. In section~4, we give a second proof that, while elementary,
makes use of continued fractions. This proof applies to an uncountable set of $\al$ including, for example,
${2}/({e-1})$. Moreover, it can be adapted to show that there are infinitely many $N$ such that
    \begin{equation*}
    \big|S_N\big(\sqrt{2}\big)\big| 
                > \frac {\log N}{2\log (1 + \sqrt 2)} + 0.78.
    \end{equation*}
This means that as $N\to\infty$ the first constant in~\eqref{eq.main} is sharp (at least for $\al=\sqrt{2}$) and
the second cannot be improved even to $3/4$.

The situation for rational $\al$ is more clear. As shown by the third author~\cite{Serbinowska.2003}, for
rational $\al$ the limit
    $$
    \lim_{N\to\infty} \frac{S_N(\al)}{N}
    $$
is well-defined. Moreover, if this limit is zero, then $S_N(\al)$ is bounded and periodic.

\begin{figure}[t]
    \begin{center}
        \begin{picture}(340,145)
            \put(332,13){$N$}
            \put(24,6){$0$}
            \put(72,6){$40$}
            \put(122,6){$80$}
            \put(170,6){$120$}
            \put(221,6){$160$}
            \put(272,6){$200$}
            \put(6,139){$S_N(\sqrt{2})$}
            \put(10,126){$3$}
            \put(10,108){$2$}
            \put(10,89){$1$}
            \put(10,71){$0$}
            \put(2,52){$-1$}
            \put(2,34){$-2$}
            \put(2,16){$-3$}
            \ifpdf
                \put(0,0){\includegraphics{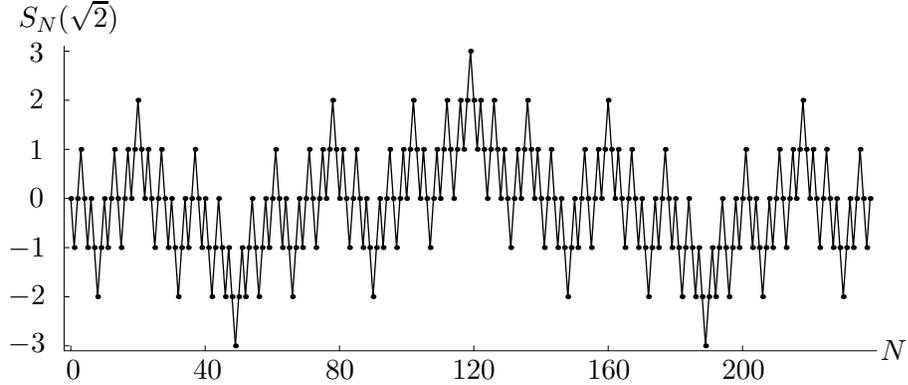}}
            \else
                \put(0,0){\includegraphics{SNroot2pic}}
            \fi
        \end{picture}
    \end{center}
    \caption{The points $\big(N,S_N(\sqrt{2})\big)$ have been connected for visual clarity.\label{figure.SNroot2pic}}

\end{figure}

In Figure~\ref{figure.SNroot2pic}, we show the points $\big(N,S_N(\sqrt{2})\big)$ when $0\le N \le 238$. It is
already apparent in this figure that $S_N(\sqrt{2})$ is not behaving like a random walk: there are never three
consecutive ``up'' steps. Also, the graph is symmetric around the peak at $N=119$ (i.e.,
$S_N(\sqrt{2})=S_{238-N}(\sqrt{2})$ for those $N$ pictured in Figure~\ref{figure.SNroot2pic}).

\begin{figure}[t]
    \begin{center}
        \begin{picture}(360,140)
            \put(310,71){$\log R_k(\sqrt{2})$}
            \put(38,55){$1$}
            \put(98,55){$4$}
            \put(159,55){$7$}
            \put(217,55){$10$}
            \put(278,55){$13$}
            \put(338,55){$16$}
            \put(18,130){$k$}
            \put(8,120){$9$}
            \put(8,102){$6$}
            \put(8,82){$3$}
            \put(8,63){$0$}
            \put(0,45){$-3$}
            \put(0,25){$-6$}
            \put(0,5){$-9$}
            \ifpdf
                \put(0,0){\includegraphics{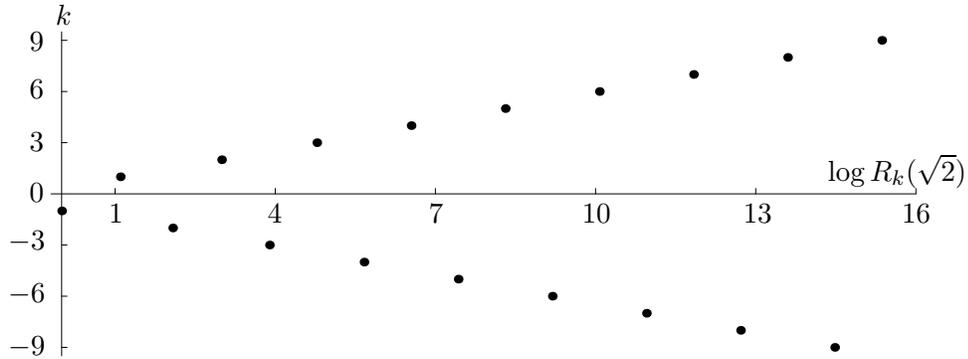}}
            \else
                \put(0,0){\includegraphics{SNroot2records}}
            \fi
        \end{picture}
    \end{center}
    \caption{The points $\big(\log R_k(\sqrt{2}),k\big)$ when $-9 \le k \le 9$, except $k=0$.\label{figure.SNroot2records}}
\end{figure}

If we restrict our attention to just the {\em record-holders}---those $N$ for which $S_N(\sqrt{2})$ takes on a
value for the first time---another aspect of the structure of $S_N(\sqrt{2})$ becomes apparent. For the sake of
rigor, we define the record-holder at the integer $k$ by
    $$
    R_k(\al) = \min\bigl\{ N\ge0 \colon S_N(\al)=k\bigr\}.
    $$
In Figure~\ref{figure.SNroot2records}, we plot the points $\big(\log R_k(\sqrt{2}),k\big)$ when $-9\le k \le 9$
(except $k=0$). The points approach two lines; we show in Corollary~\ref{cor.third.idea} that $R_k(\sqrt{2})$ is
$\floor{\frac14 (\sqrt{2}+1)^{2k+1}}$ if $k$ is positive and $\floor{\frac14 (\sqrt{2}+1)^{-2k}}$ if $k$ is
negative.

\begin{figure}[t]
    \begin{center}
        \begin{picture}(360,140)
            \put(310,91){$\log R_k(\pi)$}
            \put(38,112){$1$}
            \put(98,112){$4$}
            \put(159,112){$7$}
            \put(217,112){$10$}
            \put(278,112){$13$}
            \put(338,112){$16$}
            \put(18,130){$k$}
            \put(12,120){$4$}
            \put(10,103){$0$}
            \put(4,85){$-4$}
            \put(4,67){$-8$}
            \put(0,49){$-12$}
            \put(0,30){$-16$}
            \put(0,12){$-20$}
            \ifpdf
                \put(0,0){\includegraphics{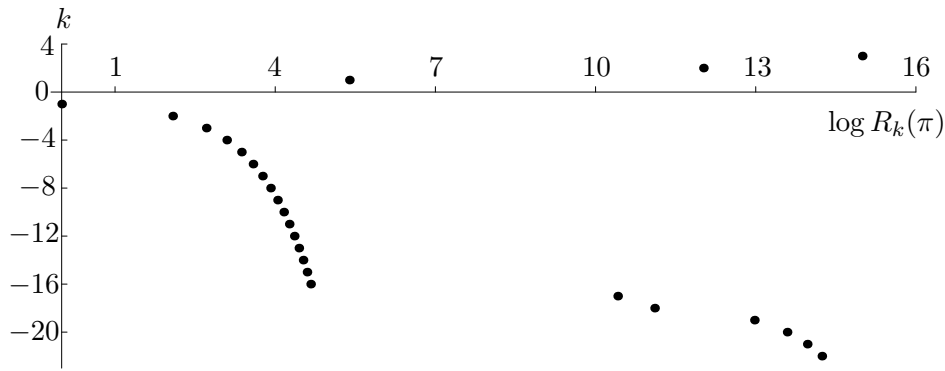}}
            \else
                \put(0,0){\includegraphics{SNpirecords}}
            \fi
        \end{picture}
    \end{center}
    \caption{The points $\big(\log R_k(\pi),k\big)$, when $-22 \le k \le 3$, except $k=0$.}
    \label{figure.SNpirecords}
\end{figure}

Our proofs do not give a logarithmic bound on $|S_N(\al)|$ for general $\alpha$; indeed, for $\al=\pi\approx
3.14159$ we do not believe that a logarithmic bound is correct. When $N\le 10^7$, computations reveal that
$-22\le S_N(\pi) \le 3$. We are unaware of any nontrivial bound on $|S_N(\pi)|$. The record holders $R_k(\pi)$
are plotted in Figure~\ref{figure.SNpirecords}. We note that the asymmetry and irregular clumping of points in
Figure~\ref{figure.SNpirecords} seems to be more typical than the orderliness depicted in
Figure~\ref{figure.SNroot2records}.

In section~2, we prove that $|S_N(\al)|\le C_\al \log N$ (for quadratic irrational $\al$) using the theory of
discrepancy, which we define but do not pursue further. Section~2 is provided for historical background and as a
hook into related literature; the remainder of this article is logically independent of section~2. In section~3,
we give an elementary bound on $|S_N(\al)|$ for a countable set of $\al$. In section~4, we introduce the
required facts and definitions about continued fractions and give a formula for $S_N(\al)$ for an uncountable
set of $\al$. From this formula, we get upper and lower bounds on the growth of $S_N(\al)$ and prove the
formulas for $R_k(\sqrt{2})$ stated earlier. We conclude in section~5 with some questions that we have not been
able to answer.

Before we begin the analysis, we introduce some notation. The natural numbers begin at 1 (i.e.,
$\N=\{1,2,3,\dots\}$). The fractional part of $x$ is given by $\fp{x}:=x-\floor{x}$. We make use of big-O
notation only with respect to $N$; that is, $f(N,\al)=\bigO{g(N)}$ if there is a constant $C$ such that
$|f(N,\al)| \le C g(N)$ for all sufficiently large $N$. Note that $C$ may depend on $\al$ but does not depend on
$N$. We also make lavish use of Iverson's notation
    $$
    \tf{Q} :=
        \left\{%
        \begin{array}{ll}
            1 & \hbox{if the statement $Q$ is true;} \\
            0 & \hbox{if the statement $Q$ is false.} \\
        \end{array}%
        \right.
    $$
See Knuth~\cite{Knuth} for an eloquent, award-winning argument in favor of this notation.\footnote{See
\url{http://www.maa.org/awards/ford.html} for a description of the Lester R. Ford Award for mathematical
exposition and the list of award-winning articles.} Iverson used $(Q)$, and Knuth uses $[Q]$, but we find
$\tf{Q}$ visually appealing, especially in light of the other frequent uses for parentheses and brackets.

\section{DISCREPANCY.}
An enduring topic in number theory has been the discrepancy shown between the expected and actual behavior of a
sequence in the interval $[0,1)$. For example, if a sequence $(w_n)$ is uniformly random, one would expect that
it is in the interval $[0,x)$ with frequency $x$. The difference between the actual and expected behavior is
formally measured by
    $$D_N^\ast\big(a_n,x\big):= \bigg| x - \frac 1N \sum_{n=1}^N \tf{w_n<x} \bigg|,$$
and the discrepancy of the sequence is denoted
    $$D_N^\ast\big( a_n \big) := \sup_{0\le x<1} D_N^\ast\big( w_n,x \big).$$
Discrepancy plays no role in the proofs of~\eqref{eq.main} given in sections~3 and~4. However, the powerful
results (stated below) of Behnke, Ostrowski and Hardy \& Littlewood, and Khintchine regarding the discrepancy of
$(\fp{n\al})$ imply bounds on $S_N(\al)$ that are more general but less precise than those given elsewhere in
this article.

We now show the connection between $S_N(\al)$ and the discrepancy of fractional part sequences. First, observe
that $\floor{n\al}$ is even exactly when there is an integer $k$ with $2k\le n\al < 2k+1$, which is the same as
$k\le n{\al}/{2} < k+1/2$. In other words,
    $$
    \tf{\floor{n\al} \text{ is even}}
        = \tf{k\le n{\al}/{2} < k+1/2}
        = \tf{\fp{ n {\al}/{2}} < 1/2}
    $$
With one eye on the definition of $D_N^\ast(\al)$, we now rewrite the definition of $S_N(\al)$ using
$(-1)^\floor{n\al} = \tf{\floor{n\al} \text{ is even}}-\tf{\floor{n\al} \text{ is odd}}$:
    \begin{align*}
    S_N(\al) &=  \sum_{n=1}^N \tf{\floor{n\al} \text{ is even}}-\tf{\floor{n\al} \text{ is odd}}. \\
    \intertext{We add this equation to the obvious}
    N &= \sum_{n=1}^N \tf{\floor{n\al} \text{ is even}}+\tf{\floor{n\al} \text{ is odd}}
    \end{align*}
to arrive at
    \begin{equation*}
    S_N(\al)    =   2 \bigg(\sum_{n=1}^N \tf{\floor{n\al} \text{ is even}} \bigg) - N
        =   -2N \; \bigg( \frac12 - \frac1N \sum_{n=1}^N \tf{\fp{n{\al}/{2}}< 1/2} \bigg).
    \end{equation*}
Thus,
    \begin{equation}\label{S-Discrepancy}
    |S_N(\al)|  =       2N\,D_N^\ast\big(\fp{n{\al}/{2}},1/2\big)
                \le     2N \, D_N^\ast\big(\fp{n{\al}/{2}}\big).
    \end{equation}

There are several accounts of the theory underlying discrepancy, most notably the colorful introductory book of
Hlawka~\cite{Hlawka} and the recent---and encyclopedic---treatise of Drmota and Tichy~\cite{Drmota.Tichy}. There
are three results that are of particular interest here.

\begin{itemize}
    \item
        Behnke \cite{Behnke1}, \cite{Behnke2}, \cite{Drmota.Tichy}*{Corollary 1.65}, \cite{Schoissengeier}
        classified        those $\al$ for which
        $$D_N^\ast \big( \fp{n\al}\big)= \bigO{\frac{\log N}{N}}.$$ In particular, this is true of all
        quadratic irrationals and is {\em not} true of $1/(e-1)$.
    \item
        Ostrowski and Hardy \& Littlewood~\cite{Drmota.Tichy}*{Theorem 1.51} showed that for any $\al$
        there are infinitely many $N$ with $$D_N^\ast\big( \fp{n\al}\big) > \frac{1}{100} \frac{\log N}{N}.$$
    \item
        Let $\psi$ be any positive increasing function on $\N$.
        Khintchine~\cite{Drmota.Tichy}*{Theorem 1.72}, \cite{Khintchine1}, \cite{Khintchine2} proved
        that $$D_N^\ast \big( \fp{n\al}\big) = \bigO{\frac{\log N}{N} \; \psi(\log N)}$$ for almost all real
        numbers $\al$ if and only if $\sum_{n=1}^\infty ({n\psi(n)})^{-1}$ converges.
\end{itemize}

The relations~\eqref{S-Discrepancy} imply that $S_N(\al)=\bigO{\log N}$ if $D_N^\ast (\fp{n{\al}/{2}}) =
\tbigO{\log N/{N}}$, but not necessarily vice versa. Since $\sqrt{2}$ and $\sqrt{5}+1$ are quadratic
irrationals, Behnke's result implies that $S_N(\sqrt{2})$ and $S_N(\sqrt{5}+1)$ have logarithmic bounds, but not
that $S_N({2}/(e-1))$ does. Schoi{\ss}engeier's work~\cite{Schoissengeier} can be used to find
    $$
    \limsup_{N\ge 1} \frac{D_N^\ast\big(\fp{n {\al}/{2}}\big)}{\log N}
    $$
for any specific quadratic irrational, but this is typically strictly larger than
    $$
    \limsup_{N\ge 1} \frac{|S_N(\al)|}{\log N}.
    $$
To be fair, we note that the discrepancy approach says that $S_N(\sqrt{3})$ has a logarithmic bound; this
{\em does not} follow from our work in the next sections.

In light of Schimdt's result, it is impossible to use a discrepancy bound to prove a bound on $S_N(\al)$ that is
sublogarithmic. As a result of Corollary~\ref{cor.third.idea}, however, for any function $\psi(N)$ with
$\psi(1)\ge1$ that increases to infinity we can find an $\al$ with $|S_N(\al)|\le\psi(N)$ for all $N$. For
example,
    $$
        S_N\Big(\frac{2}{e-1}\Big) = \bigO{\frac{\log N}{\log \log N}}.
    $$

Khintchine's result and~\eqref{S-Discrepancy} imply that for almost all real numbers $\al$
    $$
    S_N(\al) = \bigO{\log N \,\psi(\log N)}
    $$
if $\sum_{n=1}^\infty (n\psi(n))^{-1}<\infty$. In particular, set $\psi(n)=\sqrt{n}$ and observe that if
$f(N)=\tbigO{(\log N)^{3/2}}$, then for all large $N$ we have $|f(N)| \le (\log N)^2$. This verifies the rough
statement of~\eqref{eq:almostall}. It is conceivable that $S_N(\al)=\bigO{\log N}$ for almost all $\al$, but we
find this unlikely.

\section{\texorpdfstring{TWO IDEAS.}{Two Ideas}}
The sequence $S_N(\al)$ has several near-symmetries, which we formally state and prove as ``ideas.'' In this
section, we explicitly state and prove two of them, and show how they can be combined to prove
Theorem~\ref{Thm.two.ideas}. Inequality~\eqref{eq.main} is an immediate consequence of
Theorem~\ref{Thm.two.ideas}. We note that~Theorem~\ref{Thm.two.ideas} is substantively identical to the main
theorem in~\cite{Borwein.Gawronski}, although both our statement and proof are simpler.

\begin{thm}\label{Thm.two.ideas}
If $m$ is a positive integer, then
    \begin{equation*}
    \left|S_N\big(\sqrt{m^2+1}-m+1\big)\right| \le \frac{\log N}{2 \log ( \sqrt{m^2+1}+m)} + 1.
    \end{equation*}
\end{thm}

Our first idea is a combination of Beatty's theorem and a peculiar renormalization of $S_N(\al)$. We first state
Beatty's theorem, and for the sake of making this section entirely self-contained we give the ``book'' proof of
Ostrowski and Hyslop~\cite{1926.Beatty} first published in this \textsc{Monthly} in 1927. Then we give the
renormalization, and then our two ideas.

\begin{lem}[{\bf Beatty's theorem}]
If $\al>1$ is irrational and $ 1/\al +  1/\be=1$, then the sequences $(\floor{n\al})$ and $(\floor{n\be})$
partition $\N$.
\end{lem}

\begin{proof}
Let $N$ be a natural number. We will show that there is exactly one term of the two sequences ($n\al$) and
($n\be$) between $N$ and $N+1$. Since $N\al$ is not an integer, there are $\floor{N/\al}$ multiples of $\al$
less than $N$, and likewise there are $\floor{N/\be}$ multiples of $\be$. Therefore, from the two sequences
there are exactly $\floor{N/\al}+\floor{N/\be}$ terms strictly less than $N$. After writing
    $$
    N-2
    = N \left( \frac 1\al + \frac 1\be \right) - 2
    = \left(\frac N\al -1\right)+\left(\frac N\be-1\right)
    < \floor{\frac N\al}+\floor{\frac N\be}
    < \frac N\al +\frac N\be
    = N,
    $$
we see that $N-2 < \floor{N/\al}+\floor{N/\be}<N$. In other words, there are exactly
$\floor{N/\al}+\floor{N/\be} = N-1$ terms below $N$. Applying the same reasoning with $N+1$ in place of $N$, we
find that there are exactly $N$ terms below $N+1$. Thus, as claimed, there is exactly one term between $N$ and
$N+1$.
\end{proof}

And now for our peculiar renormalization. Let
    $$
    S(\al;x):= \frac14 + \sum_{0<n\al \le x} (-1)^{\floor{n\al}}
            = \frac14 + \sum_{n=1}^{\floor{x/\al}} (-1)^{\floor{n\al}}.
    $$
In terms of our earlier notation, $S_N(\al) = S(\al; N \al) - 1/4$. We have achieved two things with this
definition. First, we have introduced ``$1/4$'' at a strategic moment (with 20-20 hindsight). Second, the
parameter $x$ is naturally scaled for multiples of $\al$. This is important, since the sums $S(\al;x)$ and
$S(\be;x)$ (with $1/\al+ 1/\be=1$) are almost complementary. We invite the dubious reader to rewrite
Idea~\ref{Idea.Beatty} in terms of $S_N(\al)$.

\begin{idea}\label{Idea.Beatty}
If $\al>1$ is irrational and $1/\al+ 1/\be=1$, then $S(\al; x) + S(\be; x) = \pm 1/2.$
\end{idea}

\begin{proof}
Since $\al$ and $\be$ satisfy the hypotheses of Beatty's theorem, the set of all integers $\floor{n\al},
\floor{n\be}$---where $n\al$ and $n\be$ are in $(0,x]$---is the set $\{1,2,\ldots,r\}$, where $r$ may be either
$\floor x$ or $\floor x - 1$. Thus, we write the deepest equation in this section:
$$
S(\al;x) + S(\be;x) = \frac 24 + \sum_{j = 1}^r (-1)^j.
$$
The right-hand side is $\pm 1/2$, depending on the parity of $r$.
\end{proof}

Our second idea is really little more than the inauspicious observation that $(-1)^2=1$.

\begin{idea}\label{Idea.Simple} If $m$ is a natural number and $\be>2m$, then
$S(\be;x)=S(\be-2m;\tfrac{\be-2m}{\be} x)$.
\end{idea}

\begin{proof}
By hypothesis $\be>2m>0$, so we can multiply inequalities by $(\be-2m)/{\be}$:
    $$
        \tf{0<n\be \le x} = \tf{ 0 < n(\be-2m) \le \tfrac{\be-2m}{\be} x}.
    $$
Since
    $
        \floor{n(\be-2m)}=\floor{n\be}-2nm \equiv \floor{n\be} \pmod{2},
    $
we know that
    \begin{align*}
        S(\be;x) &=
            \frac14 + \sum_{0<n\be\le x} (-1)^{\floor{n\be}} \\
            &=
            \frac14 + \sum_{0<n (\be-2m) \le \tfrac{\be-2m}{\be} x} (-1)^{\floor{n(\be-2m)}} \\
            &=
        S(\be-2m;\tfrac{\be-2m}{\be} x).
    \end{align*}
\end{proof}

\begin{proofof}{Theorem~\ref{Thm.two.ideas}}
For the sake of being specific we work with $m=1$, setting $\al=\sqrt{2}$. Afterwards, we indicate the minor
changes needed for general $m$. Since we wish to apply Idea~\ref{Idea.Beatty}, and
$\frac{1}{\sqrt{2}}+\frac{1}{\sqrt{2}+2}=1$, we set $\be=\sqrt{2}+2$. In this notation, Idea~\ref{Idea.Beatty}
becomes
    $$
    |S(\al;x)| \le |S(\be;x)|+\frac12,
    $$
and Idea~\ref{Idea.Simple} (with $m=1$) becomes
    $$
    S(\be;x)=S\big(\al;(\sqrt{2}-1)x\big).
    $$
Set $\lambda=\sqrt{2}-1$. Combining these two ideas and iterating, we find that
    \begin{align*}
    |S(\al;x)| \le |S(\be;x)|+\frac12
        &= |S(\al; \lambda x)|+\frac12 \\
        &\le |S(\al; \lambda^2 x)| + \frac12+\frac12 \\
        & \;\;\vdots \\
        &\le |S(\al; \lambda^k x)| + \frac k2.
    \end{align*}

If $k>{\log_{\lambda}({\al}/{x})}$, then $\lambda^k x <\al$, so \mbox{$S(\al;\lambda^k x)=1/4$}. Setting
\mbox{$k=1+\floor{\log_{\lambda}({\al}/{x})}$} leads us to the inequality
    $$
    |S(\al;x)| \le \frac 14+\frac{k}{2} = \frac14+\frac 12
    \left({1+\floor{\frac{\log({\al}/{x})}{\log \lambda}}}\right)
        \le \frac{\log (x/\al) }{2 \log (\lambda^{-1})}+ \frac34.
    $$
Since $|S_N(\al)|=\left|S(\al;N\al)-1/4\right| \le |S(\al;N\al)|+ 1/4$,
    \begin{equation*}
    |S_N(\al)|
    \le \frac{\log N}{2\log (\lambda^{-1})}+1,
    \end{equation*}
as claimed. For general $m$, set $\al_m=\sqrt{m^2+1}+m-1$, $\be_m=\sqrt{m^2+1}+m+1$, and
$\lambda_m=\sqrt{m^2+1}-m$. The proof goes through verbatim.
\end{proofof}

\section{\texorpdfstring{THIRD IDEA.}{Third Idea}}
    \label{sec.Third.Idea}
Our third idea starts from the observation that if $\al$ is rational, say $\al= p/q$, then $(-1)^\floor{n\al}$
is periodic. In particular, if $N>q$ then
    \begin{align*}
    S_N\big(p/q\big)
        &= \sum_{n=1}^{q} (-1)^{\floor{np/q}} + \sum_{n=1}^{N-q} (-1)^\floor{(q+n) p/q}
        \\
        &= \sum_{n=1}^{q} (-1)^\floor{n p/q} + \sum_{n=1}^{N-q} (-1)^{p+\floor{n p/q}}
        = S_q \big(p/q\big) + (-1)^p S_{N-q}\big(p/q \big).
    \end{align*}
This allows us to replace $S_N$ with the shorter sums $S_q$ and $S_{N-q}$, and if $N-q>q$, we can replace
$S_{N-q}$ with an even shorter sum.

If $\al$ is irrational, then $(-1)^\floor{n\al}$ is not periodic, but if $\al$ is ``close'' to $p/q$, then
$(-1)^\floor{n\al}$ should be ``close'' to periodic. This is our third idea. Before we state it quantitatively,
though, we need to introduce simple continued fractions and the continued fraction expansion. Most books on
elementary number theory have chapters on continued fractions. The base-$\al$ continued fraction expansion of an
integer (which we define below) is dealt with in~\cite{MR86d:11016} and~\cites{Schoissengeier}, among other
places.

\subsection{Continued fractions.}
Throughout this article, $a_i$ always signifies an integer; if $i>0$, then $a_i$ is positive. We define the
function $[ \,\cdot\, ]$ by $[a_0]=a_0$ and for $r\ge 1$
    $$
    [a_0;a_1,a_2,\dots,a_r] =
        a_0 + \frac{1}{[a_1;a_2,\dots,a_r]} =
        a_0+\cfrac{1}{a_1+\cfrac{1}{a_2 +  \ddots   + \cfrac{1}{a_r} }}.
    $$
For every rational $p/q$ there is a unique $r$ and a unique sequence $a_0,a_1,\dots,a_r$ with $a_r\ge2$ and $p/q
= [a_0;a_1,\dots,a_r]$.

The limit
    $$
    [a_0;a_1,a_2,\dots]:=\lim_{r\to\infty} [a_0;a_1,a_2,\ldots,a_r]
    $$
always exists and is always irrational. Moreover, for each irrational $\al$ there is a unique sequence
$a_0,a_1,\dots$ with $\al=[a_0;a_1,\dots]$; we call $[a_0;a_1,\dots]$ the {\em continued fraction} of $\al$, and
the $a_t$ ($t\ge0$) the {\em partial quotients} of $\al$.

Given $\al$, one can compute the $a_i$ by noting that $a_0=\floor{\al}$ and inductively using the equation
$\al=a_0+[a_1;a_2,\dots]^{-1}$. For example, $[a_1;a_2,\dots]=(\al-a_0)^{-1}$, so $a_1=\floor{(\al-a_0)^{-1}}$,
and now $[a_2;a_3,\dots]=((\al-a_0)^{-1}-a_1)^{-1}$. The easiest concrete example is $\al=\sqrt{2}+1$. We have
$a_0=\floor{\sqrt{2}+1}=2$, and $\sqrt{2}+1=2+[a_1;a_2,\dots]^{-1}$, whence
$[a_1;a_2,\dots]=(\sqrt{2}-1)^{-1}=\sqrt{2}+1$ and $a_1=\floor{\sqrt{2}+1}=2$. Iterating, we find that
$\sqrt{2}+1=[2;2,2,2,\dots]$. We note the minor variation
    $$\sqrt{2}=[1;2,2,2,\ldots]$$
and the generalization
    $$\alpha_m=\sqrt{m^2+1}-m+1=[1;2m,2m,2m,\ldots] \quad (m\in\N),$$
which we leave as an exercise for the reader. A more difficult example is
    $$
    \frac{2}{e-1} = [1;6,10,14,18,\ldots] ,
    $$
which the very industrious can derive from the Taylor expansion of $e^x$ at $x=1$:
    $$
    e = \sum_{n=0}^\infty \frac{1}{n!}
        =  1 + \tfrac12(1+\tfrac13(1+\tfrac14(1+\cdots))).
    $$

A difficult theorem of Behnke~\cite{Drmota.Tichy}*{Corollary 1.65} states that $D_N^\ast(\al) = \tbigO{ {\log
N}/{N}}$ if and only if the sequence $\left(n^{-1} \sum_{i=1}^n a_i\right)$ is bounded, where
$\al=[a_0;a_1,a_2,\dots]$. This is true of $\alpha_m/2$, and indeed of every quadratic irrational, but is not
true of ${1}/{(e-1)}$. In fact, the set of $\al$ to which this applies has measure zero, but there are few
well-known irrationals for which it is known whether the sequence $\left(n^{-1} \sum_{i=1}^n a_i\right)$ is
bounded. For example, two famous problems are to determine the status of $\pi$ (see \cite{pi-approximation}) and
$2^{1/3}$ (see~\cite{UPINT}*{Problem F22}). We showed in section~2 that $|S_N(\al)| \le 2N \,
D_N^\ast\big(\fp{\al / 2}\big)$, so this tells us that $S_N(\al)=\bigO{\log N}$ if the partial quotients of
$\al/2$ are bounded in average. There are $\al$, however, such that $S_N(\al)$ has a logarithmic bound yet
$2N\,D_N^\ast\big(\fp{n\al/2}\big)$ does not.

We now inductively define two sequences using the partial quotients of $\al$.
    \begin{align*}
    p_{-2}&=0,  & p_{-1} &= 1,   & p_i &= a_i p_{i-1}+p_{i-2},\\
    q_{-2}&=1,  & q_{-1} &= 0,   & q_i &= a_i q_{i-1}+q_{i-2},
    \end{align*}
The $q_i$ are called the {\em continuants} of $\al$. The remarkable (albeit elementary) fact here is that
    $$
    \frac{p_i}{q_i} = [a_0;a_1,a_2,\dots,a_i].
    $$
The rationals $p_i/q_i$ are called the {\em convergents} to $\al$. The general utility of continued fractions
lies not in the convergence of the convergents to $\al$, but in the fact that $p_i/q_i$ is the closest rational
to $\al$ with denominator less than $q_{i+1}$. This is so important that we state a strong form of this
principle explicitly as Lemma~\ref{cf.is.good}. For a proof, we refer the reader to \cite{HW}*{proof of Theorem
182}, \cite{NMZ}*{Theorem 7.13}, or \cite{Rosen}*{Theorem 10.15}.

\begin{lem}\label{cf.is.good}
If $ 0<n < q_{i+1}$, then $|q_{i}\al-p_i| \le \fp{n\al} \le 1-|q_{i}\al-p_i|$, with equality only if $n= q_i$.
\end{lem}

There are several equivalent definitions of the base-$\al$ continued fraction expansion (CFE) of a
nonnegative integer $N$. Fix $\al$, and let $q_i$ denote the continuants of $\al$. We define the CFE of $N$
to be the lexicographically first sequence $(Z_i)$ of nonnegative integers satisfying
    \begin{equation}\label{def.CFE}
    N = \sum_{i=0}^\infty Z_i q_i.
    \end{equation}
In other words, write $N$ as a sum of continuants greedily (i.e., always using the largest possible), and set
$Z_i$ to be the number of times you used $q_i$. The following definitions are equivalent and useful:
\begin{itemize}
\item
    The CFE of 0 is $0,0,0,\dots$, and the CFE of $N-q_I$ is
        $$
        Z_0, \dots, Z_{I-1}, Z_I-1,Z_{I+1},\dots,
        $$
    where $q_I$ is the largest continuant less than or equal to $N$.
\item
    The $Z_i$ are nonnegative integers satisfying~\eqref{def.CFE} and
    $\sum_{i=0}^j Z_i q_i < q_{j+1}$ whenever $j\ge 0$.
\item
    The $Z_i$ are integers satisfying~\eqref{def.CFE} and
    $0\le Z_i \le a_{i+1}-\tf{Z_{i-1}>0}$ whenever $i\ge0$.
\end{itemize}

The continuants $q_0,q_1,q_2,\dots$ are positive and increasing. Therefore, if $I$ is such that $q_I>N$, then
$Z_i=0$ once $i\ge I$. We denote the CFE of $N$ by
    $$N=(Z_I,Z_{I-1},\dots,Z_1,Z_0).$$
Note that we have reversed the order of the $Z_i$; this is analogous to the custom of writing ${\bf
4}+{\bf3}\cdot 10 +{\bf2}\cdot 10^2$ as 234.

We now compute the CFEs of 100 and of $10^{11}$ with $\al=2/(e-1)=[1;6,10,14,\ldots]$. Minding our $p$s and
$q$s, we make the following table.
$$\begin{array}{|c|ccccccccc|}
\hline
 i &  0 &  1 &  2 &  3 &  4 &  5 &  6 &  7 & 8 \\
\hline
 p_i &  1 &  7 &  71 &  1001 &  18089 &  398959 &  10391023 &  312129649 & 10622799089 \\
\hline
 q_i &  1 &  6 &  61 &  860 &  15541 &  342762 &  8927353 &  268163352 & 9126481321 \\
\hline
\end{array}$$
The largest continuant not larger than 100 is 61; continuing greedily we find that $100={\bf 1} \cdot 61+
{\bf 6} \cdot 6+{\bf 3}\cdot 1$. Thus, $100=(1,6,3)$. We also have $100={\bf 16}\cdot 6+{\bf 4}\cdot1$, but
$100\not=(16,4)$ since this is not a valid CFE (the sequence $16,4,0,0,\dots$ is lexicographically after
$1,6,3,0,\dots$). A bit more arithmetic reveals that
    $$
    10^{11} = (10,32,17,6,8,15,11,9,0),
    $$
a fact that we make use of in another example below.

\begin{idea}\label{idea.Brown} If $p/q$ is the convergent to $\al$ with $q (<N)$ maximal, then
    $$
    S_N(\al) = S_q(\al)+(-1)^p S_{N-q}(\al).
    $$
\end{idea}

\begin{proof}
Define $I$ by $q=q_I$. We show that $\floor{\fp{n\al}+q\al-p}=0$ when $0<n<q_{I+1}$ by demonstrating that
    \begin{equation}\label{eq:1}
    0<\fp{n\al}+q\al-p < 1
    \end{equation}
for $n$ in this range. If $n\not=q$, then from Lemma~\ref{cf.is.good} we learn that
    $$
    -(q\al-p)\le |q\al-p| < \fp{n\al} < 1-|q\al-p| \le 1-(q\al-p),
    $$
which verifies the inequalities~\eqref{eq:1}. If $n=q$, then either $q\al-p=\fp{n\al}$ or $q\al-p =
\fp{n\al}-1$. In the first case, $0<q\al-p<1/2$ and $\fp{n\al}+q\al-p=2(q\al-p)$ lies in $(0,1)$; in the second
case, $-1/2<q\al-p<0$ and $\fp{n\al}+q\al-p=2(q\al-p)+1$ belongs to $(0,1)$. Either way, we have verified the
inequalities~\eqref{eq:1}.

Since $N-q<N \le q_{I+1}$, we have
    \begin{multline*}
    S_N(\al)-S_{q}(\al)
        =  \sum_{n=1}^{N-q} (-1)^\floor{(n+q)\al}
        =  \sum_{n=1}^{N-q} (-1)^{\floor{\floor{n\al}+\fp{n\al}+q\al-p+p}} \\
        =  \sum_{n=1}^{N-q} (-1)^\floor{n\al}(-1)^{\floor{\fp{n\al}+q\al-p}} (-1)^p
        =  (-1)^p \sum_{k=1}^{N-q} (-1)^{\floor{n\al}}
        =  (-1)^p S_{N-q}(\al).
    \end{multline*}
\end{proof}

\subsection{A formula for $S_N(\al)$.}
We are now in position to establish the following theorem. First we state it, then we give an example of it, and
then we prove it. In the next subsection, we use Theorem~\ref{Thm.third.idea} to describe the record-holders
$R_k(\al)$ for those $\al$ satisfying the hypothesis of Theorem~\ref{Thm.third.idea}.

\begin{thm}\label{Thm.third.idea} Let $\al=[a_0;a_1,a_2,\dots]$, with $a_0$ odd and
all other $a_i$ even. If a natural number $N$ has base-$\al$ CFE $(Z_I, Z_{I-1} \ldots, Z_1, Z_0)$, then
    \begin{equation*}
    S_N(\al)
        =-\sum_{i=0}^I \bigg( \prod_{t=i+1}^I (-1)^{Z_t} \bigg) \tf{i\text{ even}} \tf{Z_i \text{ odd}}.
    \end{equation*}
\end{thm}

We note that this formula can be quickly pulled out of Brown's decomposition (which we won't state) of Sturmian
words (which we won't define); indeed, this is how the authors first found Theorem~\ref{Thm.third.idea}. A
similar result, but applicable only to $\sqrt2$, appears in~\cite{MR2000b:05009}*{section 5}. The proof there
uses ideas similar to those used by Brown~\cite{MR94g:11051}. Our proof of Theorem~\ref{Thm.third.idea} borrows
from the proof of Brown's decomposition given by the first author in~\cite{OBryant.2003}.

An example will illustrate the power this result harnesses. We compute $S_N(\al)$ with $N=10^{11}$ and
    $
    \al = 2/(e-1).      
    $
We noted earlier that
    $$
        10^{11} = (10,32,17,6,8,15,11,9,0),
    $$
and now observe that $Z_i$ is odd only for $i$ in $\{1,2,3,6\}$. Thus
    $$
    \tf{i\text{ even}}\tf{Z_{i} \text{ odd}} = \tf{i \in \{2,6\}}.
    $$
We apply Theorem~\ref{Thm.third.idea} to get
    \begin{align*}
    S_{100000000000}\bigg(\frac{2}{e-1}\bigg)
        &= -\sum_{i=0}^8 \bigg( \prod_{t=i+1}^8 (-1)^{Z_t} \bigg) \tf{i\text{ even}} \tf{Z_i \text{ odd}}\\
        &= -\bigg( \prod_{t=2+1}^8 (-1)^{Z_t} + \prod_{t=6+1}^8 (-1)^{Z_t} \bigg)\\
        &= -\left( (-1)^{15+8+6+17+32+10} + (-1)^{32+10} \right)
        =-2.
    \end{align*}
A 1 Ghz PC running Mathematica would take about twenty days to compute this naively.

\begin{proofof}{Theorem~\ref{Thm.third.idea}}
Since $\al$ is fixed throughout this proof, and to avoid subscripts with subscripts, we simplify $S_N(\al)$ to
$S(N)$. By definition $p_0=a_0$, $p_1=a_0a_1+1$, and $p_i=a_i p_{i-1}+p_{i-2}$, ensuring that all $p_i$ are odd.
Idea~\ref{idea.Brown} tells us that
    \begin{equation*}
    S({q_i})
        =S({q_{i-1}})-S({q_i-q_{i-1}})=S({q_{i-1}})-\big(S({q_{i-1}})-S({q_i-2q_{i-1}})\big)
        =S({q_i-2q_{i-1}}),
    \end{equation*}
whence
    \begin{equation*}
    S(q_i)
        = S({q_i-2q_{i-1}})
        = S(q_i-4q_{i-1})
        =\dots
        = S(q_i-a_i q_{i-1})
        = S(q_{i-2})
    \end{equation*}
and so
    \begin{equation*}
    S(q_i)
        =S(q_{i-2})
        =S(q_{i-4})
        =\dots
        =
            \left\{%
            \begin{array}{ll}
                S(q_{-1}) & \hbox{if $i$ is odd,} \\
                S(q_0) & \hbox{if $i$ is even.} \\
            \end{array}%
            \right.
    \end{equation*}
Now
    $
    S({q_{-1}})= S(0) = 0
    $
and
    $$
    S(q_0)=S(1)=(-1)^\floor{\al}=(-1)^{a_0}=-1,
    $$
which completely solves the problem when $N$ is a denominator of a convergent to $\al$. Note that this gives
    $$
    S(q_i)=-\tf{i\text{ even}},
    $$
which matches with the formula stated in the theorem. This will serve both as a basis for induction, and as a
key step in the induction itself.

Suppose that the formula has been proved for all arguments less than $N=(Z_I, \ldots, Z_0)$, with $I$ chosen so
that $Z_I\not=0$. We have two cases: either $Z_I>1$ or $Z_I=1$.

If $Z_I>1$ (the easy case, so we do it first), then applying Idea~\ref{idea.Brown} twice yields (writing $q$
in place of $q_I$)
    \begin{equation*}
    S(N)=S(q)-S(N-q)=S(q)-\big(S(q)-S(N-2q)\big)=S(N-2q).
    \end{equation*}
Since the formula we are proving cares only about the parity of $Z_i$, we need only to apply the induction
hypothesis to $N-2q_I =(Z_I-2,Z_{I-1},\ldots,Z_0)$ to complete the proof.

If $Z_I=1$, then we have (again with $q=q_I$)
    \begin{align}
    S(N) &= S(q)-S\big((Z_{I-1}, \ldots ,Z_0)\big) \notag \\
         &= S(q)+
            \sum_{i=0}^{I-1}
                \bigg( \prod_{t={i+1}}^{I-1} (-1)^{Z_t} \bigg)
                \tf{i\text{ even}} \tf{Z_i \text{ odd}}. \label{eq:2}
    \end{align}
We remark that we have used the induction hypothesis applied to $(Z_{I-1}, \dots, Z_0)$, along with the fact
that the formula is not affected if we pad the left of the CFE with zeros (since it may happen that $Z_{I-1}=0$,
for example). We can expand the $S(q_I)$ term by observing that $S(q)=-\tf{I\text{ even}}$ by our base case,
$\tf{Z_I\text{ odd}}=1$ since $Z_I=1$, and $\prod_{t=I+1}^I (-1)^{Z_t}=1$ because it is an empty product. Thus,
    $$
    S(q)    = -\bigg(\prod_{t=I+1}^I (-1)^{Z_t}\bigg)\tf{I\text{ even}}\tf{Z_I\text{ odd}}.
    $$
Further, since $Z_I=1$, $\prod_{t={i+1}}^{I-1} (-1)^{Z_t}=-\prod_{t={i+1}}^{I} (-1)^{Z_t}$ when $0\le i < I$, so
    \begin{equation*}
    \sum_{i=0}^{I-1}
                \bigg( \prod_{t={i+1}}^{I-1} (-1)^{Z_t} \bigg)
                \tf{i\text{ even}} \tf{Z_i \text{ odd}}
            =
    -\sum_{i=0}^{I-1}
                \bigg( \prod_{t={i+1}}^{I} (-1)^{Z_t} \bigg)
                \tf{i\text{ even}} \tf{Z_i \text{ odd}}.
    \end{equation*}
Equation~\eqref{eq:2} now reduces to the formula as stated in the theorem.
\end{proofof}

\subsection{The record-holders.}
Recall that the record-holder at $k$ is defined by
    $$
        R_k(\al) = \min\bigl\{ N\ge0 \colon S_N(\al)=k\bigr\}.
    $$
Obviously $R_0(\al)=0$. The following corollary considers nonzero $k$ and $\al$ satisfying the hypotheses of
Theorem~\ref{Thm.third.idea}:

\begin{cor}\label{cor.third.idea}
Let $\al=[a_0;a_1,a_2,\dots]$, with $a_0$ odd and all other $a_i$ even. If $q_i$ are its continuants, then
    $$
    R_k(\al) =  q_{2k-1}\tf{k>0} + \sum_{i=0}^{2|k|-2} q_i.
    $$
\end{cor}

\begin{proof}
We briefly discuss the case in which $k<0$; the ``positive $k$'' case is similar and left to our most diligent
readers. Consider
    $$
    N_k := \sum_{i=0}^{2|k|-2} q_i.
    $$
It is not obvious that the CFE of $N_k$ is $(1,1,\dots,1)$. To verify this we must show that $\sum_{i=0}^j q_i <
q_{j+1}$ for each $j$, which requires use of the hypothesis that $a_i\ge 2$ ($i>0$).

Assuming this, we have $Z_i=1$ when $0\le i \le 2|k|-2$, and
    $$
    \prod_{t=i+1}^{2|k|-2} (-1)^{Z_t} \tf{i\text{ even}} =
        (-1)^{2|k|-2-i} \tf{i\text{ even}} =
        \tf{i\text{ even}}.
    $$
In light of this, the conclusion of Theorem~\ref{Thm.third.idea} reduces to the statement that
    $$
    S_N(\al)    =   -\sum_{i=0}^{2|k|-2} \tf{i\text{ even}}  =   k,
    $$
so $R_k(\al) \le N_k$.

Now suppose that $R_k(\al)=M=(Z_{I^\prime},Z_{{I^\prime}-1},\dots,Z_0)\le N_k$, where ${I^\prime} \le 2|k|-2$
(since $M\le N_k$). By taking the absolute value of each side of the formula given in
Theorem~\ref{Thm.third.idea} and invoking the inequality $\tf{Z_i\text{ odd}} \le 1$, we find that
    $$
    |k| =
            |S_M(\al)|
        \le \sum_{i=0}^{I^\prime} \tf{i\text{ even}}
        =   \floor{\frac {I^\prime}2}+1
        \le \frac {I^\prime}2 + 1,
    $$
and consequently that ${I^\prime}\ge 2|k|-2$. Thus, ${I^\prime}=2|k|-2$, and Theorem~\ref{Thm.third.idea} now
gives
    \begin{align*}
    k   &=
            -   \sum_{i=0}^{2|k|-2} \bigg( \prod_{t=i+1}^{2|k|-2} (-1)^{Z_t} \bigg)
                                    \tf{i\text{ even}} \tf{Z_i\text{ odd}}          \\
        &=
            -   \sum_{s=0}^{|k|-1} \bigg( \prod_{t=2s+1}^{2|k|-2} (-1)^{Z_t} \bigg)
                                    \tf{Z_{2s}\text{ odd}}.
    \end{align*}
The sum over $s$ has $|k|$ terms, each with absolute value 0 or 1, and the sum is $k$. It follows that,
    $$
    1   =   \tf{Z_{2s}\text{ odd}}  \prod_{i=2s+1}^{2|k|-2} (-1)^{Z_i}
    $$
when $0\le s \le -k-1$. For $s=|k|-1$ this implies that $Z_{2|k|-2}$ and $Z_{2|k|-1}$ are both odd, and taking
successively smaller values for $s$ informs us that all $Z_i$ are odd ($0\le i \le 2|k|-2$). Thus $M \ge
\sum_{i=0}^{2|k|-2} q_i = N_k$.
\end{proof}

If $\alpha=\sqrt{m^2+1}-m+1$, as in Theorem~\ref{Thm.two.ideas}, then we have explicit formulas for $q_i$ and
$R_k(\alpha)$. It is simply a matter of arithmetic to turn these formulas into bounds on $|S_N|$ of the type
given in Theorem~\ref{Thm.two.ideas}.

The next corollary, our finale, is in pleasant contrast to Theorem~\ref{Thm.two.ideas}.
\begin{cor} The value of $R_k(\sqrt{2})$ is given by
    $$
    R_k(\sqrt{2}) =
        \left\{%
        \begin{array}{ll}
            \floor{\tfrac14 (1+\sqrt{2})^{2k+1}} & \hbox{if $k>0$,} \\
                                                                \\
            \floor{\tfrac14 (1+\sqrt{2})^{-2k}} & \hbox{if $k<0$.} \\
        \end{array}%
        \right.
    $$
In particular, there are infinitely many $N$ with
    $$
    |S_N(\sqrt{2})| > \frac{\log (4N)}{2\log(1+\sqrt{2})}
                    > \frac{\log N}{2\log(1+\sqrt{2})} + 0.78.
    $$
\end{cor}

\begin{proof}
The continuants are given by $q_0=1$, $q_1=2$, and $q_i = 2q_{i-1}+q_{i-2}$. Solving this recurrence leads to
    $$
    q_i = \frac{(1+\sqrt 2)^{i+1}-(1-\sqrt 2)^{i+1}}{2 \sqrt 2}.
    $$
If $k<0$, we conclude that
    $$
    R_k(\sqrt2) = \sum_{i=0}^{2|k|-2} q_i = -\frac12+ \frac{(1+\sqrt{2})^{-2k}}{4} +  \frac{(1-\sqrt{2})^{-2k}}{4}
        =   \floor{\frac{(1+\sqrt{2})^{-2k}}{4}}.
    $$
We leave the positive $k$ case to the reader.

Now take $k<0$, and set $N= R_k(\sqrt2) < (1+\sqrt{2})^{2|k|}/4$. Solving for $|k|$, we find that
    $$
    |S_N(\sqrt2)| = |k| > \frac{\log (4N)}{2\log(1+\sqrt{2})}.
    $$
\end{proof}

\section{\texorpdfstring{FOUR HARDER QUESTIONS.}{Four Harder Questions}}
    \label{sec.Harder.Questions}
We now state four problems that we have been unable to resolve.
\begin{itemize}
    \item
    What happens for other $\al$? For example, we know from section~2 that if $\al$ is
    any quadratic irrational, then $S_N(\al)=\bigO{\log N}$. Is this the correct type of growth for all
    quadratic irrationals? What are the necessary and sufficient conditions on $\al$ for $S_N(\al) =
    \bigO{\log N}$ in terms of the continued fraction expansion of $\al$? It seems unlikely that $S_N(\al)$ is
    $\bigO{\log N}$ for almost all $\al$, but a proof of this is elusive.

    \item
    What can be said for a given base $m$ about the number of $\floor{n\al}$ ($1\le n \le N$)
    that belong to the various
    congruence classes modulo $m$? There are $m-1$ free parameters here, because the sum is always $N$.
    For $m=2$, this is a problem equivalent to the one that we have worked on in this article; for $m\ge3$ it is
    different. As a starting point, one may wish to consider the sum
        $$\sum_{n=1}^N (e^{2\pi i /m})^{\floor{n\al}}.$$

    \item
    Our ideas do not generalize in any really straightforward manner to handle nonhomogeneous sequences. For
     example, is
        $$
        \left|\sum_{n=1}^N (-1)^{\floor{n\sqrt{2}+1/2}}\right| \leq \frac{\log N}{2\log(1+\sqrt{2})}+\frac32
        $$
    when $N>2$?

    \item
    The problem of bounding $S_N(\sqrt{2})$ originally arose~\cite{Ruderman} in studying the convergence of
        $$
        \sum_{n=1}^\infty \frac{(-1)^{\floor{n\sqrt 2}}}{n}.
        $$
     Is it possible to evaluate this infinite sum? Is it
     rational? Schmuland~\cite{Schmuland} studies the distribution of $\sum_{n=1}^\infty w_n/n$, where the
     $w_n$ are independent, identically distributed random
     variables taking the values 1 and $-1$ with equal probability. Does the sum $\sum_{n=1}^\infty
     (-1)^{\floor{n/\alpha}}/n$, where $\alpha$ is chosen uniformly from $(0,1)$, have the same distribution?
\end{itemize}

\vspace{12pt} \small \noindent {\bf ACKNOWLEDGMENT.}\, The first author thanks the National Science Foundation
(grant DMS-0202460) for its support. The authors thank Dennis Eichhorn for carefully proofreading the
manuscript.

\vspace{12pt} \noindent {REFERENCES}

\vspace{12pt}

\noindent {\bf KEVIN O'BRYANT} is a University of Illinois at Urbana-Champaign alumnus, and after a stint as an
NSF postdoc (DMS 0202460) at the University of California, San Diego, he has taken a permanent position at the
City University of New York. His research is mostly in combinatorial number theory and is strongly motivated by
the philosophy of experimental mathematics. His hobbies include backgammon, nonstandard analysis, and drinking
coffee.
\newline {\em Department of Mathematics, City University of New York, CSI, Staten Island, NY 10314
\newline kevin@member.ams.org}

\vspace{12pt} \noindent {\bf BRUCE REZNICK} got his degrees from Caltech (B.S., 1973) and Stanford (Ph.D.,
1976), and his been at the University of Illinois at Urbana-Champaign since 1979. He is interested in
combinatorial problems in number theory, algebra, analysis, geometry, often involving polynomials. He still
roots for the Chicago Cubs, but has forgotten why.\newline
    {\em
    Department of Mathematics, University of Illinois at Urbana-Champaign, 1409 W. Green St., Urbana, IL
    61801\newline
    reznick@math.uiuc.edu
    }

\vspace{12pt} \noindent {\bf MONIKA SERBINOWSKA} enjoyed growing up in  Warsaw, Poland. She obtained her M.S.
from Warsaw University in 1990 under the supervision of Wojciech Guzi\`{n}ski. In 1996 she obtained a Ph.D. from
the University of Utah for work in theoretical statistics under the direction of Lajos Horv\'{a}th.  Currently, she
is an assistant professor at Weber State University. She is now pursuing questions from number theory, continued
fractions and geometry. She has three sons and she is an avid bridge player, hiker and skier.
    \newline
    {\em
    Department of Mathematics, Weber State University, Ogden, UT 84408 \newline
    mserbinowska@weber.edu
    }

\end{document}